\documentclass[12pt]{amsart}
\advance\textheight\topskip
\usepackage[OT2,T1]{fontenc}
\newcommand{\mifody}{%
  \renewcommand\rmdefault{wncyr}%
  \renewcommand\sfdefault{wncyss}%
  \renewcommand\encodingdefault{OT2}%
  \normalfont
  \selectfont}
\usepackage{times}
\allowdisplaybreaks

\usepackage[matrix,arrow,curve,dvips]{xy} 

\usepackage{afterpage}
\usepackage{graphicx}
\newcounter{mapsfigure}

\usepackage[dvips, bookmarks=false, hyperfigures=false,
setpagesize=false]{hyperref}

\usepackage[mathscr]{eucal}
\usepackage{mathabx}
\usepackage{mathptmx}

\newcommand{\pone}{(p,1)}

\newcommand{\State}{\mathscr{F}}
\newcommand{\Vertex}[2]{{#1}^{\{#2\}}_{\vphantom{h}}}

\newcommand{\Bbin}[2]{{\Shuffle^{#1}_{#2}}}
\newcommand{\Bfac}[1]{\mathfrak{S}_{#1}}

\newcommand{\Shuffle}{\mathop{\text{\mifody\sf Sh}}\nolimits}

\newcommand{\algO}{\mathscr{O}_{p,1}}

\newcommand{\Co}[1]{{}#1{}\tensor{}}
\newcommand{\pbw}[3]{F_1^{#1} F_3^{#2} F_2^{#3}}
\newcommand{\PBW}[3]{\textsc{b}(#1,#2,#3)}
\newcommand{\Pbw}[3]{\textsc{b}(#1,#2,#3)}

\newcommand{\tensor}{\mathbin{\otimes}}

\newcommand{\phalpha}{\varphi_{\alpha}}
\newcommand{\phbeta}{\varphi_{\beta}}
\newcommand{\Ea}{\mathscr{E}_{\alpha}}
\newcommand{\Eb}{\mathscr{E}_{\beta}}

\newcommand{\ttimes}{\pmb{\times}}
\newcommand{\ttimescirc}{{\circ}\kern-4pt{\ttimes}}

\newcommand{\hw}[1]{{\bigl|}#1\bigr\rangle}
\newcommand{\q}{\mathfrak{q}}

\newcommand{\adj}{%
  \mathchoice{\mathbin{\blacktriangleright}}%
  {\mathbin{\mbox{\small${\blacktriangleright}$}}}%
  {\mathbin{{\blacktriangleright}}}%
  {\mathbin{{\blacktriangleright}}}}

\newcommand{\Aint}[1]{\langle#1\rangle}
\newcommand{\Afac}[1]{\langle#1\rangle!\,}
\newcommand{\Abin}[2]{\mathchoice%
  {\Abinm{#1}{#2}\,}{\Abinmm{#1}{#2}\,}%
  {\Abinmm{#1}{#2}}{\Abinmm{#1}{#2}}}
\newcommand{\Abinm}[2]{\mbox{\footnotesize$\displaystyle
    \genfrac{\langle}{\rangle}{0pt}{}{#1}{#2}$}}
\newcommand{\Abinmm}[2]{\genfrac{\langle}{\rangle}{0pt}{}{#1}{#2}}

\newcommand{\Dynkin}[3]{\xymatrix@C28pt@1{\circ\ar@{-}[r]^(.55){#2}\ar@{}^{#1}[]&\circ\ar@{}^{#3}[]}}

\newcommand{\bref}[1]{\textbf{\textup{\ref{#1}}}}

\newcommand{\mfrac}[2]{\raisebox{.8pt}{\mbox{\small$\displaystyle\frac{#1}{#2}$}}}
\newcommand{\ffrac}[2]{\raisebox{.5pt}{\mbox{\footnotesize$\displaystyle\frac{#1}{#2}$}}}
\newcommand{\fffrac}[2]{\raisebox{.9pt}{\mbox{\scriptsize$\displaystyle\frac{#1}{#2}$}}}
\newcommand{\half}{%
  \mathchoice{\ffrac{1}{2}}{\frac{1}{2}}{\frac{1}{2}}{\frac{1}{2}}}

\newcommand{\algW}{\mathcal{W}}
\newcommand{\W}{\mathscr{W}}
\newcommand{\Wa}{\mathscr{W}_{\alpha}}
\newcommand{\Wb}{\mathscr{W}_{\beta}}
\newcommand{\Wba}{\mathscr{W}_{\beta\alpha}}
\newcommand{\Wab}{\mathscr{W}_{\alpha\beta}}
\newcommand{\Waba}{\mathscr{W}_{\alpha\beta\alpha}}
\newcommand{\Wbab}{\mathscr{W}_{\beta\alpha\beta}}
\newcommand{\Waabb}{\mathscr{W}_{\alpha\alpha\beta\beta}}

\newcommand{\Fa}{F_{\alpha}}
\newcommand{\Fb}{F_{\beta}}

\newcommand{\HHyd}{{}\mbox{\small${}^{H}_{H}$}\YDname}
\newcommand{\YDname}{\mathscr{Y\kern-1ptD}}

\newcommand{\oC}{\mathbb{C}}
\newcommand{\oZ}{\mathbb{Z}}

\newcommand{\Nich}{\mathfrak{B}}
\newcommand{\BX}{\Nich(X)}

\numberwithin{equation}{section}
\makeatletter
\@addtoreset{equation}{section}
\@addtoreset{subsubsection}{section}

\def\@secnumfont{\bfseries}
\def\subsubsection{\@startsection{subsubsection}{3}%
  \z@{.5\linespacing\@plus.7\linespacing}{-.5em}%
  {\normalfont\bfseries}}
\def\paragraph{\@startsection{paragraph}{4}%
  \z@\z@{-\fontdimen2\font}%
  \normalfont\bfseries}
\def\subparagraph{\@startsection{subparagraph}{5}%
  \z@\z@{-\fontdimen2\font}%
  \normalfont\bfseries}

\makeatother

\swapnumbers

\theoremstyle{definition}

\newtheorem{rem}[subsubsection]{Remark}

\begin{document}

\title{A note on the ``logarithmic-$\algW_3$'' octuplet algebra and
  its Nichols algebra}

\author[Semikhatov]{AM Semikhatov}

\address{Lebedev Physics Institute, Moscow 119991, Russia}

\begin{abstract}
  We describe a Nichols-algebra-motivated construction of an octuplet
  chiral algebra that is a ``$\algW_3$-counterpart'' of the triplet
  algebra of $(p,1)$ logarithmic models of two-dimensional conformal
  field theory.
\end{abstract}

\maketitle
\thispagestyle{empty}

\section{Introduction}
Logarithmic models of two-dimensional conformal field theory can be
defined as centralizers of Nichols
algebras~\cite{[STbr],[c-charge]}. \ For this, the generators $F_i$ of
a given Nichols algebra $\Nich(X)$ with diagonal
braiding~\cite{[Nich],[Wor],[Lu-intro],[Rosso-inv],[AG],[AS-pointed],
  [Andr-remarks],[AS-onthe],[Heck-Weyl],[Heck-class],[AHS],[ARS],
  [Ag-0804-standard],[Ag-1008-presentation],[Ag-1104-diagonal]} are to
be realized as
\begin{equation*}
  F_i = \oint e^{\alpha_i . \varphi},\quad 1\leq i \leq
  \text{rank}\equiv \theta,
\end{equation*}
where $\varphi(z)$ is a $\theta$-plet of scalar fields and
$\alpha_i\in\oC^{\theta}$ are chosen so as to reproduce the given
braiding coefficients $q_{i,j}$ in
\begin{equation*}
  \Psi: F_i\tensor F_j\mapsto q_{i,j} F_j\tensor F_i,\quad
  1\leq i,j\leq\theta.
\end{equation*}
The coefficients are standardly arranged into a braiding matrix
$(q_{i,j})_{\substack{1\leq i\leq\text{rank}\\
    1\leq j\leq\text{rank}}}$.  The relation between the braiding
matrix and the screening momenta is postulated~\cite{[c-charge]} in
the form of equations
\begin{equation*}
  q_{j,j} = e^{i\pi \alpha_j.\alpha_j},\qquad
  q_{j,k}q_{k,j} = e^{2i\pi \alpha_j.\alpha_k}
\end{equation*}
and the logical-``or'' conditions
\begin{equation*}
  a_{i,j} \alpha_i.\alpha_i = 2\alpha_i.\alpha_j,
  \quad\mbox{\small$\bigvee$}\quad
 (1-a_{i,j})\alpha_i.\alpha_i = 2
\end{equation*}
imposed for each pair $i\neq j$ and involving the Cartan matrix
$a_{i,j}$ associated with the given braiding matrix (see, e.g.,
\cite{[HS]} and the references therein).

In this note, we describe some details related to the construction of
the octuplet algebra~\cite{[c-charge]} that can be considered a
``logarithmic extension'' of the $\algW_3$ algebra~\cite{[Z]}
similarly to how the triplet algebra~\cite{[Kausch],[GK+],[FHST]} is a
``logarithmic extension'' of the Virasoro algebra.  The starting point
is a particular item in Heckenberger's list of rank-2 Nichols algebras
with diagonal braiding (which is item \textbf{5.7}(1)
in~\cite{[Hel]})---the braiding matrix
\begin{equation}\label{eq:qij}
  q_{ij}=\begin{pmatrix}
    \q^2 & \q^{-1} \\
    \q^{-1} & \q^2
  \end{pmatrix},  
\end{equation}
where $\q^2$ is a primitive $2p$th root of unity.  We choose
\begin{equation}\label{the-q}
  \q = e^{\frac{i\pi}{p}}
\end{equation}
with $p=2,3,\dots$.  This choice leads to $\pone$-type logarithmic
CFT{}
models~\cite{[Kausch],[GK+],[FHST],[FGST],[FGST2],[AM-3],[AM-latt]},
in contrast to $(p,p')$ models that follow if $\q$ is chosen as
$e^{\frac{i\pi p'}{p}}$ instead.  The main expectation associated with
$\pone$-type models is that their representation categories are
``\textit{very} closely related''~\cite{[FGST2],[NT],[TW]} to an
appropriate representation category on the algebraic side, which in
the braided case is some category of Yetter--Drinfeld
$\Nich(X)$-modules (cf.~\cite{[b-fusion]}).  In this paper, we
therefore proceed along two routes: (i)~describing the structure of
the $\Nich(X)$ algebra associated with~\eqref{eq:qij} (solely with the
choice in~\eqref{the-q}) and its suitable Yetter--Drinfeld modules,
and (ii)~discussing some properties of the octuplet algebra that
centralizes this~$\Nich(X)$.  None of the two directions is pursued to
the point where they actually meet (which would mean constructing a
functor), but the results presented here hopefully bring us somewhat
closer to that point.

\section{The Nichols algebra}
\subsection{Presentation for $\Nich(X)$}
We first recall the presentation of the relevant Nichols algebra, as a
quotient of the tensor algebra.  Our starting point is a
two-dimensional braided vector space $X$ with the preferred basis
$F_1$, $F_2$ and the above braiding matrix in this basis.  The Nichols
algebra $\BX$ is the quotient by a graded
ideal~$\mathscr{I}$~\cite{[Ag-1008-presentation],[Hel]},
\begin{equation}\label{Nich-W3}
  \BX=T(X)/ \bigl([F_1,[F_1,F_2]],\ [F_2,[F_2,F_1]],\ F_1^{p},\
  [F_2,F_1]^{p},\ F_2^{p}\bigr),\quad \dim\Nich(X)=p^3,
\end{equation}
If $p=2$, the double-bracket generators of the ideal are absent.  The
brackets here denote $\q$-commutators determined by the braiding
matrix: $[F_1,F_2]= F_1 F_2 - \q^{-1}F_2 F_1$, $[F_2,F_1]= F_2 F_1 -
\q^{-1}F_1 F_2$, and so on by multiplicativity of the ``$q$''-factor,
whence the two double commutators in the ideal are explicitly given by
\begin{align*}
  [F_1,[F_1,F_2]]&=F_1^2 F_2 -(\q+\q^{-1})F_1 F_2 F_1 + F_2 F_1^2,\\
  [F_2,[F_2,F_1]]&=F_2^2 F_1 -(\q+\q^{-1})F_2 F_1 F_2 + F_1 F_2^2.
\end{align*}

A PBW basis in $\Nich(X)$ is given by $\pbw{r}{t}{s}$, $0\leq
r,s,t\leq p-1$ \cite{[Hel]}, where
\begin{equation*}
  F_3=[F_2,F_1].
\end{equation*}
The double-bracket relations in the ideal can also be rewritten as
$F_2 F_3 = \q F_3 F_2$ and
$F_3 F_1 = \q F_1 F_3$.

Multiplication in $\Nich(X)=T(X)/\mathscr{I}$ is the one induced by
``concatenation'' in the tensor algebra, $X^{\otimes m}\tensor
X^{\otimes n}\to X^{\otimes(m+n)}$,
$(x_1,\dots,x_m)\tensor(y_1,\dots,y_n)\mapsto
(x_1,\dots,x_m,y_1,\dots,y_n)$.  It is then relatively straightforward
to show that the multiplication table of the PBW basis elements is
\begin{multline}\label{the-product}
  (\pbw{r_1}{t_1}{s_1})(\pbw{r_2}{t_2}{s_2}) ={}\\
  \sum_{i=0}^{\min(s_1,r_2)}
  \q^{t_1 (r_2 - i) + t_2 (s_1 - i) - s_1 r_2 + i (1 + i)/2}
  \Afac{i}\Abin{s_1}{i}\Abin{r_2}{i}
  \pbw{r_1 + r_2 - i}{t_1 + t_2 + i}{s_1 + s_2 - i}. 
\end{multline}
Comultiplication is by ``deshuffling,'' determined by the defining
property of a braided Hopf algebra and the fact that $F_1$ and $F_2$
are primitive.

\subsection{$\Nich(X)$ as a subalgebra in $T(X)$}
For any Nichols algebra $\Nich(X)$, the graded ideal $\mathscr{I}$
such that $\Nich(X)=T(X)/\mathscr{I}$ is known to be the kernel of the
total braided symmetrizer map in each grade, $\Bfac{n}:X^{\otimes
  n}\to X^{\otimes n}$.  Mapping by $\Bfac{n}$ in each grade therefore
results in an equivalent description of $\Nich(X)$ with multiplication
given by the shuffle product
\begin{equation*}
  {*}{}:{}(x_1,\dots,x_m)\tensor(y_1,\dots,y_n)\mapsto
  \Bbin{}{m,n}(x_1,\dots,x_m,y_1,\dots,y_n),
\end{equation*}
and comultiplication by deconcatenation (see~\cite{[STbr]} for the
definition of shuffles and the braided symmetrizer; the only
notational difference is that $*$ is not used for the shuffle product
there).

We let $\PBW{r}{t}{s}$ be the image of $\pbw{r}{t}{s}$ under the map
by the braided symmetrizer, or more precisely,
\begin{equation*}
  \PBW{r}{t}{s}
  =\ffrac{1}{\Afac{r} \Afac{s} \Afac{t} (1 - \q^2)^{t}}\Bfac{r+2t+s}
  (\pbw{r}{t}{s}).
\end{equation*}
In particular,
\begin{alignat*}{2}
  \PBW{1}{0}{0}&=F_1,&\qquad  \PBW{2}{0}{0}&= F_1 F_1,\\
  \PBW{0}{0}{1}&=F_2,&  \PBW{1}{0}{1}&=F_1 F_2 + \q^{-1} F_2 F_1,\\
  &&\PBW{0}{0}{2}&= F_2 F_2,\\
  &&\PBW{0}{1}{0}&=-\q^{-2} F_2 F_1.
\end{alignat*}

\subsubsection{}\label{sec:product}
The shuffle product of $\PBW{r_1}{t_1}{s_1}$ and $\PBW{r_2}{t_2}{s_2}$
follows from~\eqref{the-product}:
\begin{multline}\label{the-product-Sh}
  \PBW{r_1}{t_1}{s_1}*\PBW{r_2}{t_2}{s_2} ={}\\
  \sum_{i=0}^{\min(s_1,r_2)}\!\!\!
  \Abin{r_1 + r_2 - i}{r_1} \Abin{s_1 + s_2 - i}{s_2}
  \mfrac{
    (1 - \q^2)^i \Afac{t_1 + t_2 + i}}{\Afac{t_1}\Afac{t_2}\Afac{i}}\;
  \q^{t_1 (r_2 - i) + t_2 (s_1 - i) - s_1 r_2 + i (i + 1)/2}\\*
  {}\times
  \PBW{r_1 + r_2 - i}{t_1 + t_2 + i}{s_1 + s_2 - i}.
\end{multline}
and the coproduct is
\begin{align}\label{Delta}
  \Delta:\PBW{r}{t}{s}&{}\mapsto \sum_{j=0}^{r} \sum_{m=0}^{s}
  \sum_{k=0}^{t} \sum_{i=0}^{k} (-1)^i \q^{-i (i + 3)/2 + (k - m - 2
    i) j + m (t - i - k)}
  \\
  & \notag \quad{}\times\Abin{i + j}{i} \Abin{i + m}{i} \Afac{i}
  \Co{\PBW{r - j}{k - i}{i + m}} \PBW{j + i}{t - k}{s - m},
  \\
  \intertext{where terms with the lowest grades in the first tensor
    factor are} \notag &=1\tensor\PBW{r}{t}{s} + F_1\tensor
  \PBW{r-1}{t}{s}
  \\
  \notag &\quad{}+\q^{t - r} F_2\tensor \PBW{r}{t}{s-1} - \q^{-r - 2}
  \Aint{r+1} F_2\tensor \PBW{r + 1}{t - 1}{s}
  \\
  \notag &\quad{}+\dots
\end{align}
(the dots stand for terms $\PBW{r'}{t'}{s'}\tensor\PBW{r''}{t''}{s''}$
with $r'+2t'+s'\geq2$).

\begin{rem}
  Although this is obvious, we note explicitly that the ``Serre
  relations''---the double $q$-commutators in the ideal---are resolved
  in terms of the shuffle product in the sense that the relations
  \begin{align*}
    F_1 * F_1 * F_2 -(\q+\q^{-1})F_1 * F_2 * F_1 + F_2 * F_1 * F_1=0,\\
    F_2 * F_2 * F_1 -(\q+\q^{-1})F_2 * F_1 * F_2 + F_1 * F_2 * F_2=0
  \end{align*}
  hold identically for the shuffle product defined by the braiding
  matrix~\eqref{eq:qij}.
\end{rem}

\subsubsection{}\label{sec:antipode}
The action of the antipode on the PBW basis elements is defined by the
formulas
\begin{align*}
  S(\PBW{r}{0}{0}) &= (-1)^r \q^{r (r - 1)} \PBW{r}{0}{0},\\
  S(\PBW{0}{t}{0}) &= \sum_{i=0}^{t} (-1)^t
  \q^{\half i (i - 1) - (i + 3) t + t^2} \Afac{i}\, \PBW{i}{t - i}{i},\\
  S(\PBW{0}{0}{s}) &= (-1)^s \q^{s (s - 1)} \PBW{0}{0}{s}
\end{align*}
and by the fact that $S$ is a braided antiautomorphism:
\begin{equation*}
  S(\PBW{r}{t}{s})=
  \q^{r t - r s + t s}
  S(\PBW{0}{0}{s}) *
  S(\PBW{0}{t}{0}) * S(\PBW{r}{0}{0}).
\end{equation*}

\subsection{Vertex operators and Yetter--Drinfeld $\Nich(X)$ modules}\label{sec:YDmodules}
Multivertex $\Nich(X)$ module comodules, which are Yetter--Drinfeld
modules, were defined in~\cite{[STbr]}.  We here realize simple
Yetter--Drinfeld modules of our $\Nich(X)$ in terms of one-vertex
modules.

\subsubsection{The $Y$ spaces}\label{Vn1n2}
Let $Y^{\{n_1, n_2\}}$ be a one-dimensional vector space with basis
$V^{\{n_1, n_2\}}$ and braiding $\psi:\Nich(X)\tensor Y^{\{n_1, n_2\}}
\to Y^{\{n_1, n_2\}}\tensor\Nich(X)$ and $Y^{\{n_1,
  n_2\}}\tensor\Nich(X)\to \Nich(X)\tensor Y^{\{n_1, n_2\}}$ defined
by
\begin{align*}
  \psi(F_i\tensor V^{\{n_1, n_2\}})
  &= \q^{1-n_i}\,V^{\{n_1, n_2\}}\tensor F_i, 
  \\
  \psi(V^{\{n_1, n_2\}}\tensor F_i)
  &= \q^{1-n_i} F_i\tensor V^{\{n_1, n_2\}},
\end{align*}
$i=1,2$.  Every space $\Nich(X)\tensor V^{\{n_1^1, n_2^1\}}\tensor
\Nich(X)\tensor V^{\{n_1^2, n_2^2\}}\tensor\dots\tensor V^{\{n_1^N,
  n_2^N\}}$ is a Yetter--Drinfeld $\Nich(X)$ module.  Taking the
$a_i^j$ to be generic leads to continuum families of such modules,
leaving us with no chance of a nice correspondence with any type of
``reasonably rational'' CFT{} model.  The choice of the possible
$a_i^j$ values is governed by the requirement that all of them (and
the braided vector space $X$ itself) be objects of a suitable $\HHyd$
category of Yetter--Drinfeld modules over a nonbraided Hopf
algebra~$H$.  In the case of diagonal braiding, more specifically,
$H=k\Gamma$ for an Abelian group~$\Gamma$, which can then be
considered the origin of the appropriate discreteness in the $a_i^j$
values.  We do not pursue this line in this paper, and simply assume
that the $a_i^j$ take integer values.

We consider one-vertex modules $\Nich(X)\tensor V^{\{n_1, n_2\}}$ and
for brevity write
\begin{align*}
  \Vertex{\PBW{r}{t}{s}}{n_1, n_2}&=\PBW{r}{t}{s}\tensor
  V^{\{n_1, n_2\}}\in\Nich(X)\tensor Y^{\{n_1, n_2\}},
  \\
  \intertext{and, in particular,}
  \Vertex{F_i}{n_1, n_2}&= F_i\tensor
  V^{\{n_1, n_2\}}\in\Nich(X)\tensor Y^{\{n_1, n_2\}}
\end{align*}
(but $\Vertex{\PBW{0}{0}{0}}{n_1, n_2} = 1\tensor V^{\{n_1, n_2\}}$ is
normally written as $V^{\{n_1, n_2\}}$).

\subsubsection{Left adjoint action}
The formulas for the product, coproduct, and antipode in
\bref{sec:product}--\bref{sec:antipode} allow calculating the left
adjoint action of the $\Nich(X)$ generators on one-vertex modules:
\begin{multline*}
  F_1 \adj \Vertex{\Pbw{r}{t}{s}}{n_1, n_2}
  =
  \Aint{r+1}
  (1 - \q^{2(r - s + t + 1 - n_1)}) 
  \Vertex{\Pbw{r + 1}{t}{s}}{n_1, n_2}\\*
  {}- \q^{2 r - 2 s + t - 2 n_1 + 3} \Aint{t + 1} 
  (1 - \q^2) \Vertex{\Pbw{r}{t + 1}{s - 1}}{n_1, n_2}
\end{multline*}
and
\begin{multline*}
  F_2 \adj \Vertex{\Pbw{r}{t}{s}}{n_1, n_2}
  =
  \q^{1 - r} \Aint{t + 1} (1 - \q^2)
  \Vertex{\Pbw{r - 1}{t + 1}{s}}{n_1, n_2}\\*
  {}+ 
  \q^{t - r} \Aint{s+1} (1 - \q^{2 (s + 1 - n_2)})
  \Vertex{\Pbw{r}{t}{s + 1}}{n_1, n_2}.
\end{multline*}
These formulas depend on $n_1$ and $n_2$ only through
$(a_i\;\text{mod}\;p)$.  The $\Nich(X)$ coaction is given by literally
applying formula~\eqref{Delta} to $\Pbw{r}{t}{s}\tensor V^{\{n_1,
  n_2\}}$ (and is entirely independent of~$a_i$).

\subsubsection{Simple Yetter--Drinfeld modules}
A simple Yetter--Drinfeld $\Nich(X)$-module $\mathscr{Y}_{n_1,n_2}$ is
generated from $V^{\{n_1,n_2\}}$ under the action of $\Nich(X)$; its
dimension is given by
\begin{equation*}
  d(p,n_1,n_2)=
  \begin{cases}
    d(\overline{n_1},\overline{n_2})
    ,&
    \overline{n_1} + \overline{n_2} \leq p,
    \\
    d(\overline{n_1},\overline{n_2})
    - d(p-\overline{n_1},p-\overline{n_2})
    ,&
    \overline{n_1} + \overline{n_2}
    \geq p + 1,
  \end{cases}
\end{equation*}
where $d(n_1,n_2)=\half\,n_1\,n_2\,(n_1 + n_2)$ 
and
\begin{equation*}
  \overline{x}=
  \begin{cases}
    p,&(x\;\text{mod}\;p)=0,\\
    x\;\text{mod}\;p,&\text{otherwise}.
  \end{cases}
\end{equation*}



\section{The octuplet algebra centralizing $\Nich(X)$}%
\setcounter{figure}{1}%
We next discuss a CFT{} construction related to our $\Nich(X)$.

\subsection{Screenings and their zero-momentum centralizer}
We identify the $\Nich(X)$ generators with two screenings
\begin{equation}\label{the-scr}
  \Fa=F_1=\oint e^{\phalpha},\qquad
  \Fb=F_2=\oint e^{\phbeta},
\end{equation}
where $\phalpha(z)$ and $\phbeta(z)$ are two scalar fields whose OPEs
are defined in accordance with the braiding matrix as follows:
\begin{alignat*}{2}
  \phalpha(z)\,\phalpha(w)&=\ffrac{2}{p}\log(z-w),\quad
  &\phalpha(z)\phbeta(w)&=-\ffrac{1}{p}\log(z-w),
  \\
  &&\phbeta(z)\phbeta(w)&=\ffrac{2}{p}\log(z-w).
\end{alignat*}

It follows from the formulas in~\cite{[c-charge]} that the centralizer
(``kernel'') of screenings~\eqref{the-scr} contains a Virasoro algebra
with the central charge
\begin{equation}\label{W3-cc}
  c= 50 - \ffrac{24}{p} - 24 p
  = -\ffrac{2 (3 p - 4) (4 p - 3)}{p}.
\end{equation}
This Virasoro algebra is represented by the energy--momentum tensor
\begin{equation*}
  T(z) =
  \ffrac{p}{3} 
  \partial\phalpha\partial\phalpha(z)
  +\ffrac{p}{3} \partial\phalpha\partial\phbeta(z)
  +\ffrac{p}{3} \partial\phbeta\partial\phbeta(z)
  -(p\!-\!1) \partial^2\phalpha(z)
  -(p\!-\!1) \partial^2\phbeta(z).
\end{equation*}
In addition to the Virasoro algebra, the kernel of the screenings
contains the dimension-3 Virasoro primary field (omitting the
conventional $(z)$ arguments of fields)
\begin{multline}\label{W3-gen}
  W(z) =
  \partial\phalpha \partial\phalpha \partial\phalpha
  +
  \ffrac{3}{2}\partial\phalpha \partial\phalpha
  \partial\phbeta 
  -\ffrac{3}{2} \partial\phalpha \partial\phbeta
  \partial\phbeta
  -\partial\phbeta \partial\phbeta \partial\phbeta
  \\
  -\ffrac{9 (p - 1)}{2
    p}\partial^2\phalpha \partial\phalpha -\ffrac{9 (p -
    1)}{4 p} \partial^2\phalpha \partial\phbeta + \ffrac{9
    (p - 1)}{4 p}\partial^2\phbeta \partial\phalpha +
  \ffrac{9 (p - 1)}{2 p} \partial^2\phbeta \partial\phbeta
  \\
  + \ffrac{9 (p - 1)^2}{4 p^2} \partial^3\phalpha -\ffrac{9 (p -
    1)^2}{4 p^2} \partial^3\phbeta.
\end{multline}
The operator product of this field with itself is given by
\begin{align*}
  W(z)\,W(w) &=
  \ffrac{81  (3 p-5) (3 p-4) (4 p-3) (5 p-3)}{4 p^5}\ffrac{1}{(z-w)^6}
  - \ffrac{243}{4 p^4}\,\ffrac{(3 p - 5) (5 p - 3) T(w)}{(z-w)^4}
  \\
  &\quad{}
  - \ffrac{243}{8 p^4}\,\ffrac{(3 p - 5) (5 p - 3)\partial T(w)}{(z-w)^3} 
  +\ffrac{243}{16 p^4}\,\ffrac{8 p TT(w)
    - 9 (p - 1)^2 \partial^2 T(w)}{(z-w)^2}
  \\
  &\quad{}
  + \ffrac{243}{8 p^4}\,\ffrac{4 p (\partial T)T(w)
    - (p-1)^2 \partial^3 T(w)}{z-w},
\end{align*}
where $TT(w)$ is the normal-ordered product $T(w)T(w)$ (and similarly
for $(\partial T)T(w)$).  This OPE defines the $\algW_3$
algebra~\cite{[Z]} (also see~\cite{[BS]}).

In an equivalent description, the $\algW_3$ algebra relations for the
modes introduced as $T(z)=\sum_{n\in\oZ} L_n z^{-n-2}$ and
$W(z)=\sum_{n\in\oZ} W_n z^{-n-3}$ are
\begin{align*}
  [L_{m}, L_{n}] &= (m - n) L_{m + n} + \ffrac{1}{12}
  (50 - \ffrac{24}{p} - 24 p) (m - 1) m (m + 1)\delta_{m + n, 0},\\
  [L_{m}, W_{n}] &= (2 m - n) W_{m + n},\\
  [W_{m}, W_{n}] &= -\ffrac{81 (3 p\!-\!5) (5 p\!-\!3)}{8 p^4} (m - n)
  \Bigl( \ffrac{(m + n + 3) (m + n + 2)}{5} - \ffrac{(m + 2) (n +
    2)}{2} \Bigr)
  L_{m + n} \\
  &\quad{}+ \ffrac{243}{4 p^3} (m - n) \Lambda_{m + n}
  + \ffrac{27 (3 p\!-\!5) (3 p\!-\!4) (4 p\!-\!3) (5 p\!-\!3)}{
    160 p^5} m (m^2\!- 1) (m^2\!- 4)\delta_{m + n, 0},
\end{align*}
where 
\begin{equation*}
  \Lambda_m = \sum_{n\leq -2} L_{n} L_{m - n}
  + \sum_{n\geq -1} L_{m - n}  L_{n}
  - \ffrac{3}{10} (m + 3) (m + 2) L_{m}.
\end{equation*}

\subsection{Long screenings}
The $\algW_3$ algebra is also centralized by two ``long'' screenings
\begin{equation}\label{long-scr}
  \Ea
  =\oint e^{-p\phalpha}
  \quad\text{and}\quad
  \Eb
  =\oint e^{-p\phbeta}.
\end{equation}
Because
\begin{equation*}
  [F_i,\mathscr{E}_j]=0,
\end{equation*}
the long screenings act on the kernel of the $\Fa$ and $\Fb$, and are
therefore a useful tool in studying that kernel.

\subsection{Remark}
We note that, generally, given the screenings $F_i=\oint e^{\varphi_i}
=\oint e^{\alpha_i\cdot\varphi}$, $i=1,\dots,\theta$, the Virasoro
dimension of a vertex $e^{\mu.\varphi(z)}$ with $\mu=\sum_{i=1}^\theta
c_i\alpha_i$ is
\begin{equation*}
  \Delta(c)
  =\sum_{i=1}^\theta c_i\bigl(1-\ffrac{\alpha_i.\alpha_i}{2}\bigr)
  +\half \sum_{i,j=1}^\theta c_i c_j\,\alpha_i.\alpha_j.
\end{equation*}
We list the generators of the ideal in~\eqref{Nich-W3} together with
the vertex operators that naively (by momentum counting) correspond to
them, and with the Virasoro dimensions of these vertices:
\begin{equation}\label{corresp-dims}
  \begin{alignedat}{5}
     &[F_1,[F_1,F_2]], &\quad& [F_2,[F_2,F_1]], &\quad& F_1^{p},
     &\quad& [F_1,F_2]^{p}, &\quad& F_2^{p},
     \\
     &e^{2\phalpha(z)+\phbeta(z)},
     &\quad&e^{\phalpha(z)+2\phbeta(z)},
     &\quad&e^{p\phalpha(z)},
     &\quad&e^{p\phalpha(z)+p\phbeta(z)},
     &\quad&e^{p\phbeta(z)},
     \\
     &3,&\quad& 3,&\quad& 2p-1,&\quad& 3p-2,&\quad& 2p-1.
  \end{alignedat}
\end{equation}

\subsection{The octuplet algebra}
The field
\begin{equation*}
  \W(z)=e^{p\phalpha(z) + p\phbeta(z)},
\end{equation*}
which is the top-dimension field in~\eqref{corresp-dims}, is in the
kernel of $\Fa$ and $\Fb$ and is a $\algW_3$-primary field of
dimension $\Delta = 3p-2$ and the $W_0$ eigenvalue zero.  To describe
how it is mapped by the long screenings, we need a reminder on
$\algW_3$ singular vectors.

\subsubsection{Singular vectors in $\algW_3$ Verma
  modules}\label{xy-sing-vect}
We recall from~\cite{[BW]} (also see~\cite{[BS]} and the references
therein) that highest-weight vectors of the $\algW_3$ algebra can be
conveniently parameterized by $(x,y)$ such that
\begin{align*}
  L_m\hw{x, y} &=0,\quad m\geq 1,\\
  W_m\hw{x, y} &=0,\quad m\geq 1,\\
  L_{0} \hw{x, y} &= \Bigl(\ffrac{x^2 + y^2 + x y}{3} - \ffrac{(p -
    1)^2}{p}\Bigr)
  \hw{x, y},\\
  W_{0} \hw{x, y} &= \ffrac{1}{2 p^{3/2}} (x - y) (2 x + y) (x + 2 y)
  \hw{x, y}.
\end{align*}
The two numbers $x$ and $y$ are defined not uniquely but up to a Weyl
transformation; the Weyl group orbit of $(x,y)$ also contains $(-x,
x+y)$, $(x+y,-y)$, $(y, -x-y)$, $(-x-y,x)$, and $(-y,-x)$. \ We write
$\mathscr{V}(z)\doteq\hw{x,y}$ for any field$/$state $\mathscr{V}(z)$
that satisfies the above conditions.

In what follows, we use the conditions for the existence of singular
vectors in Verma modules of the $\algW_3$
algebra~\cite{[Mchi],[W],[BW]}.  Whenever a state can be represented
as $\hw{x,y}$ with $x=a\sqrt{p} - \fffrac{c}{\sqrt{p}}$ for integer
$a$ and $c$ such that $ac>0$, there is a singular vector on the level
$ac$ built on that state.  The singular vector has the highest-weight
parameters $(x',y')=(x-2a\sqrt{p}, y + a\sqrt{p})$. \ Similarly, if
$y=b\sqrt{p} - \fffrac{d}{\sqrt{p}}$ with $bd>0$, then a singular
vector occurs on the level~$bd$ and has the highest-weight parameters
$(x'',y'')=(x + b\sqrt{p}, y - 2b\sqrt{p})$.

\subsubsection{}
It follows that
\begin{align*}
  \W(z)=e^{p\phalpha(z) + p\phbeta(z)} &\doteq
  \hw{2\sqrt{p}-\fffrac{1}{\sqrt{p}}, 2 \sqrt{p}-\fffrac{1}{\sqrt{p}}},
\end{align*}
and hence the corresponding Verma-module state has two singular
vectors at level 2.  Both of them vanish in our free-field
realization.  Of the two fields $\Ea\W(z)$ and $\Eb\W(z)$, we
concentrate on the second; it lands in the module generated from 
\begin{align*}
  e^{p\phalpha(z)}
  &\doteq \hw{3\sqrt{p}-\fffrac{1}{\sqrt{p}},-\fffrac{1}{\sqrt{p}}}.
\end{align*}
The corresponding highest-weight state in the Verma module has
singular vectors at levels $3$ and $p-1$.  The first of these vanishes
in the free-boson realization, but the second does not, yielding just
the field $\Wb(z)=\Eb\W(z)$, as we show in
Fig.~\thefigure.\setcounter{mapsfigure}{\thefigure} We note that
\afterpage{%
  \rotatebox[origin=lB]{90}{\centering%
    \parbox{.98\textheight}{%
      \begin{gather*}
        {}\\[-2\baselineskip]
        \xymatrix@C6pt@R12pt{ &&&&&&&&&&&&&&&&&&&&&&&&&&&&&&&&&
          *{{\makebox[0pt]{$e^{-p\phalpha-p\phbeta}$\qquad\ \
              }}}&*{\circ} \ar@{-}[dddddddddddddddllllll]
          \ar@{-}[dddddddddrrrr] \ar^(.8){p-1}[ddddl] \ar[ddddr] &&
          \\
          &&&&&&&&&&&&&&&&&&&&&&&&&&&&&&&&&&&&
          \\
          &&&&&&&&&&&&&&&&&&&&&&&&&&&&&&&&&&&&
          \\
          &&&&&&&&&&&&& &*{1\!\!\!}&*{\circ} \ar@{-}[ddddddddddddllll]
          \ar@{-}[dddddddddddrrrr] \ar[ddddr]_(.8){4\!\!}
          \ar^(.85){\!\!\!2p-1}[ddddddl] \ar;[dddddd]+<12pt,4pt>
          &&&&&&&&&&&&&&&&&&&&&
          \\
          &&&&&&&&&&&&&&&&&&&&&&&
          &*{\makebox[0pt]{$e^{-p\phbeta}$}}&*{\circ}
          \ar@{-}[dddddddddddllll]\ar@{-}[ddddddddrrr]
          \ar|(.7){\;2}[dd]+<-4pt,4pt> \ar;[dd]+<4pt,4pt>
          \ar^(.8){2p-2}[dddddl] \ar[dddddr]
          \ar@/_6pt/|{\Ea}[rrrrrrrr] &&&&&&&&*{\bullet} \ar[dddddl]
          \ar[dddddrrr] \ar[dd]\ar[ddrr] &&*{\bullet}\ar[dddddr]
          \ar[dd]\ar[ddll] \ar^(.7){2p-2}[dddddlll] &&
          \\
          &&&&&&&&&&&&&&&&&&&&&&&&&&&&&\ar@{{.}{.}{>}}@/_6pt/|(.3){\Eb}[urrrrrr]&&&&&&&
          \\
          &&&&&&&&&&&&&&&&&&&&&&&&&*{{\ttimes}{\ttimes}}&&&&&&&&*{\ttimes}&&*{\ttimes}&
          \\
          &&&&&&&&&&&&&&&&*{\ttimes}&&&&&&&&&&&&\ar@{{.}{.}{>}}@/_4pt/|(.25){\Eb}[ddrrrrrrrr]&&&&&&&&
          \\
          &&&&&&&&&&&&&&&&&&&&&&&&&&&&&&&&&&&&
          \\
          &&&&&&*{\makebox[0pt]{$e^{p\phalpha}$}} &*{\circ}
          \ar@{-}[ddddddlll] \ar@{-}[dddddrr] \ar[ddddl]^(.7){p-1}
          \ar[dd]^(.7){\!\!3} \ar@/_6pt/|{\Ea}[rrrrrrr]
          &&&&&&&*{\ttimescirc}\ar[ddddl]^(.7){p-1}
          \ar@/_10pt/|{\Eb}[]+<6pt,-3pt>;[rrrrrrrrrr] &
          &*{{\kern-12pt\ttimescirc}}
          \ar@/^10pt/|{\Eb}[]+<3pt,3pt>;[rrrrrrrrrr]
          \ar^(.7){\!\!p-1}[]+<-4pt,-4pt>;[ddddr]
          &&&&&&&&*{\ttimescirc} \ar[ddddr]^(.5){\!\!p-1}
          \ar@/_20pt/|{\Ea}[rrrrrrrrrrrr] &&*{\ttimescirc} \ar[ddddl]
          \ar@/_6pt/|{\Ea}[rrrrrr] &&&&&
          &*{\ttimescirc}\ar^(.7){\!\!p-1}[];[ddddrr]
          &&&&*{\ttimescirc}\ar[];[ddddll]&&&
          \\
          &&&&&&&&&&&&&&&&&&&&&&&&&&&&&&&&&&&&
          \\
          &&&&&&&*{\ttimes}&&&\ar@{{.}{.}{>}}@/_10pt/|(.75){\Eb}[uurrrrrr]+<-10pt,-3pt>&&&&&&&&&&&&&&&&&&\ar@{{.}{.}{>}}@/_10pt/|(.7){\Eb}[uurrrrrrrr]&&&&&&&&
          \\
          &&&&&&&&&&&&&&&&&&&&&&&&&&&&&&&&&&&&
          \\
          &*{\W}&*{\bullet}\ar@{-}[ddl]\ar@{-}[ddr]
          \ar@/^6pt/|{\Eb}[rrrr]&&&&*{\bullet}\ar@/_6pt/|(.6){\Ea}[rrrrrrr]&*{\!\!\!{}^{\textstyle\Wb}}&&&&&
          &*{\bullet}\ar@/_16pt/|(.6){\Eb}[rrrrrrrrrrrr]&*{\makebox[0pt]{$\quad\Wab$}}&
          &*{\Wba\!\!\!}&*{\bullet} \ar@/_6pt/|(.45){\Eb}[rrrrrrrr]
          &&&&&&&
          &*{\bullet}\ar@/_16pt/|{\Ea}[];[rrrrrrrrr]&*{\makebox[0pt]{$\quad\Wbab$}}
          &&&&&&& &*{\bullet}&*{\makebox[0pt]{$\quad\ \;\Waabb$}}&
          \\
          &&&&&&&&&&&&&&&&&&&&&&&&&&&&&&&&&&&&\\
          &&&&&&&&&&\ar@{{.}{.}{>}}@/_10pt/|{\Eb}[uurrrrrrr]&&&&&&&&&&&&&&&&&\ar@{{.}{.}{>}}@/_10pt/|(.6){\Eb}[uurrrrrrr]&&&&&&&&&
        } \end{gather*}%
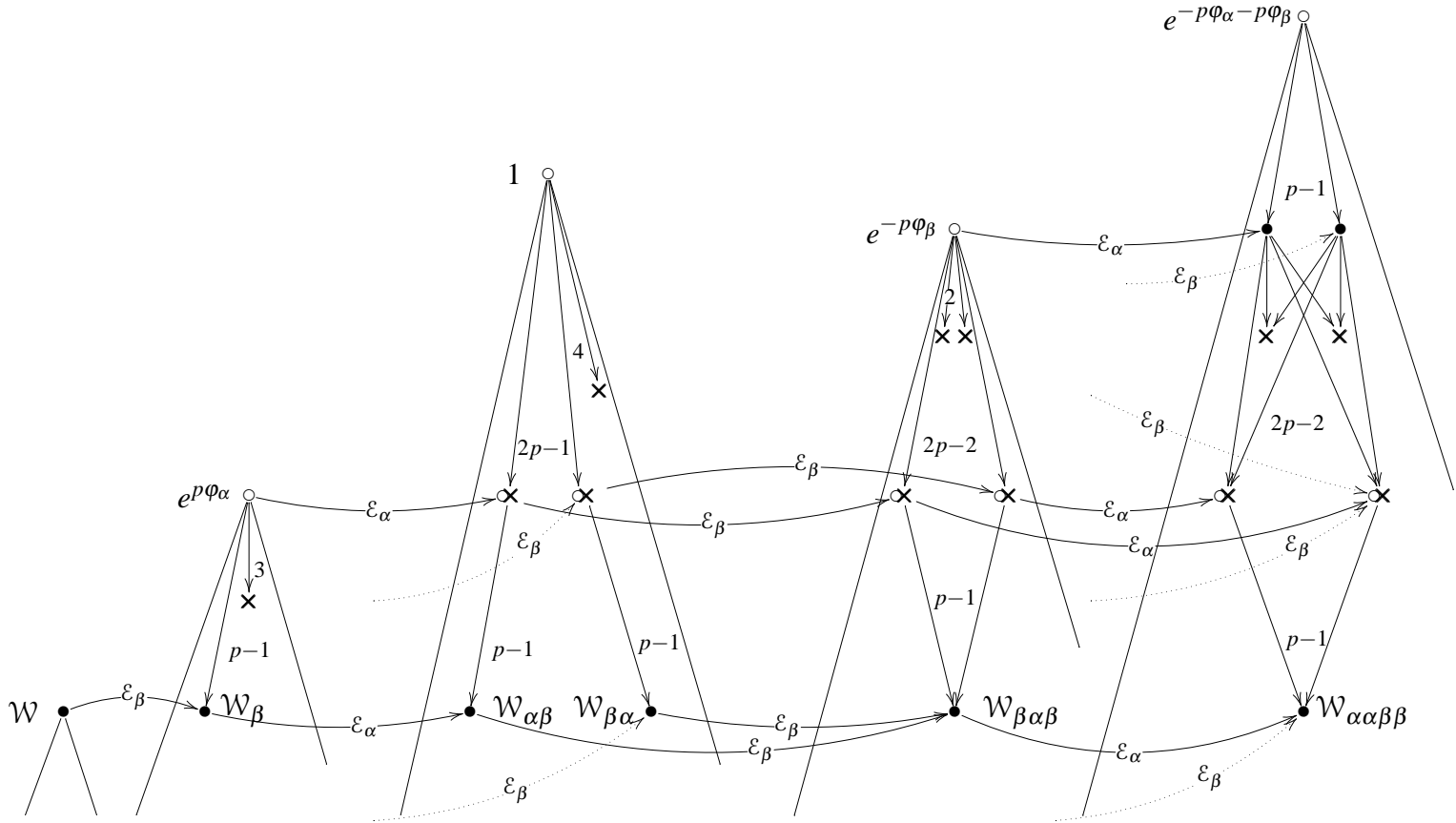
      \textsc{Figure}~\thefigure. \ \small Maps by the long screenings
      $\Ea$ and $\Eb$.  Crosses and downward arrows leading to them
      show $\algW_3$ singular vectors that vanish in the free-boson
      realization.  Bullets (and downward arrows) show nonvanishing
      states in the same grades; the relative levels of singular
      vectors are indicated at the arrows.  An open circle
      superimposed with a cross shows a vanishing $\algW_3$ singular
      vector and a ($\algW_3$-primary) state in the same grade, but
      not in the same $\algW_3$-module (and downward arrows drawn from
      such $\ttimescirc$ show singular vectors built on those primary
      states).  Two more modules---those with $e^{p\phbeta}$ and
      $e^{-p\phalpha}$ at the top---are not shown here; their
      structure repeats that of the ``$e^{p\phalpha}$'' and
      ``$e^{-p\phbeta}$'' modules with $\alpha\leftrightarrow\beta$.
      Dotted arrows show the maps by $\Ea$ and $\Eb$ \textit{from} the
      missing modules.
    }}%
  \addtocounter{figure}{1}\mbox{}\\
  \mbox{}\\
  \mbox{}}
\begin{align*}
  \W_{\beta}(z)&=\mathscr{P}_{\beta}^{[p -1]}(\partial\varphi(z))
  e^{p\phalpha(z)}
\end{align*}
with a differential polynomials in $\partial\phalpha(z)$,
$\partial\phbeta(z)$ in front of the exponential; here and hereafter,
we indicate the degree $d$ of a differential polynomial as
$\mathscr{P}^{[d]}$.

Totally similarly,
\begin{align*}
  \W_{\alpha}(z)=\Ea\W(z)&=\mathscr{P}_{\alpha}^{[p -1]}(\partial\varphi(z))
  e^{p\phbeta(z)}
  \\
  \intertext{is a descendant of}
  e^{p\phbeta(z)}
  &\doteq \hw{-\fffrac{1}{\sqrt{p}},3
    \sqrt{p}-\fffrac{1}{\sqrt{p}}}.
\end{align*}

The maps of $\Wa(z)$ by $\Eb$ and of $\Wb(z)$ by $\Ea$ are
differential polynomials (not involving exponentials).  They are not
descendants of the unit operator, however.  We have $1 \doteq
\hw{\sqrt{p}-\frac{1}{\sqrt{p}},\sqrt{p}-\frac{1}{\sqrt{p}}}$, which
implies singular vectors at levels $1$, $1$, $4$, $2p-1$, and $2p-1$.
All of these vanish in the free-field realization.  In each of the
grades where a level-$(2p-1)$ singular vector vanishes, another state
is produced as $\Ea(e^{p\phalpha(z)})$ and $\Eb(e^{p\phbeta(z)})$. \
This is shown in Fig.~\themapsfigure\ with the
$\;{\times}\kern-11.5pt\circ\ $ symbols (``a state superimposed with a
vanishing singular vector'').  Next, $\Ea(e^{p\phalpha(z)})$ and
$\Eb(e^{p\phbeta(z)})$ have singular-vector descendants on the
relative level $p-1$, which are of course the respective images of
$\Wb(z)$ and $\Wa(z)$ under $\Ea$ and $\Eb$,
\begin{equation*}
  \W_{\alpha\beta}(z)=
  \Ea\W_{\beta}(z)=\mathscr{P}_{\alpha\beta}^{[3 p - 2]}(\partial\varphi(z)),
  \qquad
  \W_{\beta\alpha}(z)=
  \Eb\W_{\alpha}(z)=\mathscr{P}_{\beta\alpha}^{[3 p -2]}(\partial\varphi(z)).
\end{equation*}

Further maps by the long screenings do not produce
$\algW_3$-descendants of the corresponding exponentials either.  We
consider $\Eb\W_{\alpha\beta}(z)$ and $\Eb\W_{\beta\alpha}(z)$.  In
the module associated with
\begin{align*}
  e^{-p\phbeta(z)}
  &\doteq
  \hw{2\sqrt{p}-\fffrac{1}{\sqrt{p}},-\sqrt{p}-\fffrac{1}{\sqrt{p}}},
\end{align*}
two singular vectors at level $2$ and two at level $2p-2$ vanish;
located at the grades of the last two are $\Eb\Ea e^{p\phalpha(z)}$
(the maps shown in Fig.\;\themapsfigure) and $\Eb\Eb
e^{p\phbeta(z)}$.\footnote{We illustrate the use
  of~\bref{xy-sing-vect}.  In the Verma module with the highest-weight
  vector $\hw{x,y}
  =\hw{2\sqrt{p}-\frac{1}{\sqrt{p}},-\sqrt{p}-\frac{1}{\sqrt{p}}}$
  associated with $e^{-p\phbeta(z)}$, one of the level-$(2p-2)$
  singular vectors exists due to the representation
  $y=\frac{p-1}{\sqrt{p}} - 2\sqrt{p}$, and therefore the singular
  vector has the highest-weight parameters
  $(x'',y'')=(-\frac{1}{\sqrt{p}},3\sqrt{p}-\frac{1}{\sqrt{p}})$,
  i.e., those of $e^{p\phbeta(z)}$.  The other level-$(2p-2)$ singular
  vector is seen immediately if we Weyl-reflect the highest-weight
  parameters to $(\tilde x,\tilde y)=(-x,x+y)$.  We then have $\tilde
  y=\frac{2(p-1)}{\sqrt{p}} - \sqrt{p}$, and hence the singular vector
  has the parameters
  $(-3\sqrt{p}+\frac{1}{\sqrt{p}},3\sqrt{p}-\frac{2}{\sqrt{p}})$.
  After the same Weyl reflection, the parameters
  $(3\sqrt{p}-\frac{1}{\sqrt{p}},-\frac{1}{\sqrt{p}})$ correspond to
  $e^{p\phalpha(z)}$.}  Now, $\Eb\Ea e^{p\phalpha(z)}$ and $\Eb\Eb
e^{p\phbeta(z)}$ have a level-$(p-1)$ singular-vector descendant each.
In our free-field realization, these two singular vectors evaluate the
same up to a nonzero overall factor, thus producing a
$\algW_3$-primary field
\begin{align*}
  \W_{\beta\alpha\beta}(z)=
  \Eb\W_{\alpha\beta}(z)&=\mathscr{P}_{\beta\alpha\beta}^{[3 p - 3]}
  (\partial\varphi(z))
  e^{-p\phbeta(z)}.
  \\
  \intertext{Everything with the replacement
    $\alpha\leftrightarrow\beta$ applies to the field}
  \W_{\alpha\beta\alpha}(z)=
  \Ea\W_{\beta\alpha}(z)&=\mathscr{P}_{\alpha\beta\alpha}^{[3 p -
    3]}(\partial\varphi(z)) e^{-p\phalpha(z)}.
\end{align*}

Finally, mapping by the long screenings once again gives a field
\begin{align*}
  \W_{\alpha\alpha\beta\beta}(z)=
  \Ea\W_{\beta\alpha\beta}(z)&=\mathscr{P}_{\alpha\alpha\beta\beta}^{[4 p - 4]}
  (\partial\varphi(z)) e^{-p\phalpha(z) - p\phbeta(z)}
\end{align*}
(which is also $\Eb\W_{\alpha\beta\alpha}(z)$ up to a factor), which
is not in the module associated with $e^{-p\phalpha(z)-p\phbeta(z)}$,
however.  In the Verma module associated with the highest-weight
vector
\begin{align*}
  e^{-p\phalpha(z) - p\phbeta(z)}
  &\doteq \hw{-\fffrac{1}{\sqrt{p}},-\fffrac{1}{\sqrt{p}}},
\end{align*}
there are two singular vectors at level $p-1$, both of which are
nonvanishing in the free-field realization and are in fact the images
of $e^{-p\phbeta(z)}$ (and $e^{-p\phalpha(z)}$; see
Fig.~\themapsfigure).  Each of these singular vectors therefore has
two level-$(2p-2)$ singular vectors, which are in fact the same pair
of singular vectors.  These two next-generation singular vectors
vanish in our free-field realization, but the maps by $\Ea$ (and by
$\Eb$) land in the same grades.  The two vectors in the image of the
long screenings share a singular-vector descendant at the
level-$(p-1)$ and this descendant is the $\Waabb(z)$ field.

\subsubsection{}We summarize the octuplet structure of $\algW_3$ primary
fields generated by long screenings from $\W(z)$:
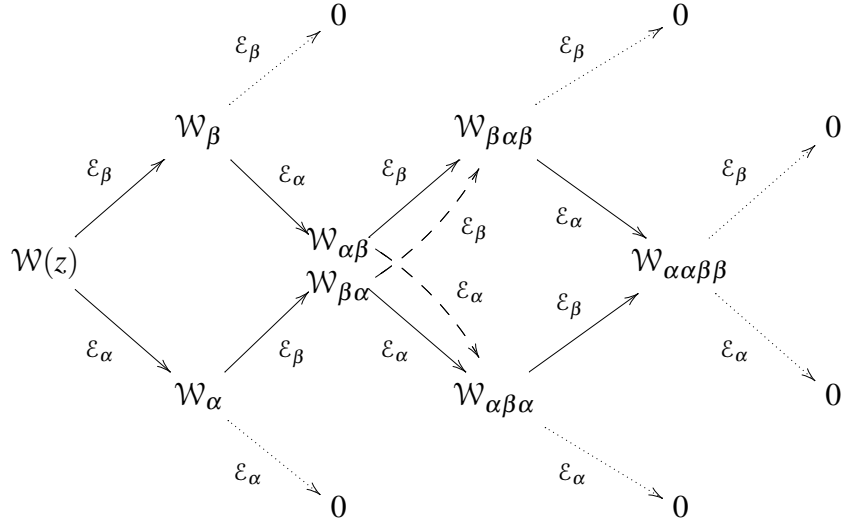
\begin{equation*}
  \xymatrix@C28pt{
    &&0&&0&\\
    &\Wb\ar@{{.}{.}{>}}^{\Eb}[ur]\ar[rd]^{\Ea}&&\Wbab\ar[rd]_{\Ea}\ar@{{.}{.}{>}}^{\Eb}[ur]&&0\\
    \W(z)\ar[ru]^{\Eb}\ar[rd]_{\Ea}&&
    *{\mbox{}\kern-6pt\begin{array}{c}\Wab\\
        \Wba
      \end{array}\kern-6pt}\ar_(.65){\Eb}@{{-}{--}{>}}@/_3pt/[]+<16pt,-4pt>;[ur]+<-6pt,-15pt>
    \ar^(.6){\Ea}@{{-}{--}{>}}@/^3pt/[]+<16pt,4pt>;[dr]+<-6pt,15pt>
    \ar[ru]^{\Eb}\ar[rd]_{\Ea}&&\Waabb\ar@{{.}{.}{>}}^{\Eb}[ur]\ar@{{.}{.}{>}}_{\Ea}[dr]&\\
    &\Wa\ar@{{.}{.}{>}}_{\Ea}[dr]\ar[ru]_{\Eb}&&\Waba\ar[ru]^{\Eb}\ar@{{.}{.}{>}}_{\Ea}[dr]&&0\\
    &&0&&0&
    }
\end{equation*}
The dashed arrows represent maps to the target field up to a nonzero
overall factor.  All the fields in the diagram are
$\algW_3$-primaries, with the same Virasoro dimension $3p-2$.

We follow~\cite{[c-charge]} in proposing these fields as generators of
the octuplet algebra $\algO$---the extended algebra of logarithmic
$\algW_3$ models.

\subsubsection{}Calculations with particular examples show the OPE
\begin{align*}
  \W(z)\,\Waabb(w)&=\mfrac{c_1\cdot 1}{(z-w)^{6p-4}} +
  \mfrac{c_2\,T(w)}{(z-w)^{6p-6}} + \mfrac{c_2/2\,\partial
    T(w)}{(z-w)^{6p-7}} + \dots
  \\
  \intertext{with \textit{nonzero} $p$-dependent coefficients (and no
    dimension-$3$ $W(w)$ field), and}
  \Wa(z)\,\Wbab(w)&=\mfrac{(-1)^{p+1} c_1\cdot 1}{(z-w)^{6p-4}} +
  \mfrac{(-1)^{p+1} c_2\,T(w)}{(z-w)^{6p-6}}
  + \mfrac{(-1)^{p+1} c_2/2\,\partial T(w)}{(z-w)^{6p-7}}+\dots,\\
  \Wb(z)\,\Waba(w)&=\mfrac{(-1)^{p+1} c_1\cdot 1}{(z-w)^{6p-4}} +
  \mfrac{(-1)^{p+1} c_2\,T(w)}{(z-w)^{6p-6}} + \mfrac{(-1)^{p+1}
    c_2/2\,\partial T(w)}{(z-w)^{6p-7}}+\dots.
\end{align*}
with \textit{nonzero} coefficients,
and the OPEs $\Wa(z)\,\Wbab(w)$ and $\Wb(z)\,\Waba(w)$ that start very
similarly.  The adjoint-$s\ell(3)$ nature of the octuplet manifests
itself in the OPEs such as
\begin{align*}
  \Wa(z)\,\Wb(w)&=\mfrac{c_3 \W(w)}{(z-w)^{3p-2}}+\dots,\\
  \Wa(z)\,\Waba(w)&=\mathcal{O}(z-w),\\
  \Wb(z)\,\Wbab(w)&=\mathcal{O}(z-w),\\
  \Waba(z)\,\Wbab(w)&=\mfrac{c'_3 \Waabb(w)}{(z-w)^{3p-2}}+\dots.
\end{align*}

\subsubsection{Some octuplet algebra representations}
To construct CFT{} counterparts of the modules introduced
in~\bref{sec:YDmodules}, we first define the ``fundamental weights''
$\omega_i$ such that $\omega_i\cdot\alpha_j=\delta_{i,j}$:
\begin{equation*}
  \omega_1=\ffrac{p}{3}(2\alpha_1 + \alpha_2),\qquad
  \omega_2=\ffrac{p}{3}(\alpha_1 + 2\alpha_2).
\end{equation*}
We let $\omega_{\alpha}(z)$ and $\omega_{\beta}(z)$ denote the
corresponding fields:
\begin{equation*}
  \omega_{\alpha}(z)=\ffrac{p}{3}(2\phalpha(z) + \phbeta(z)),\qquad
  \omega_{\beta}(z)=\ffrac{p}{3}(\phalpha(z) + 2\phbeta(z)).
\end{equation*}
Then the field
\begin{equation*}
  \State_{n_1,n_2}(z)=
  e^{\frac{1-n_1}{p}\omega_{\alpha}(z) + \frac{1-n_2}{p}\omega_{\beta}(z)}
\end{equation*}
has the same braiding with $F_i$ as $V^{\{n_1, n_2\}}$ has
in~\bref{Vn1n2}.  The dimension of $\State_{n_1,n_2}(z)$ is
\begin{equation*}
  \Delta_{n_1,n_2}=
  p - n_1 - n_2 + \frac{n_1^2+n_1 n_2+n_2^2}{3 p} - \ffrac{(p - 1)^2}{p}
\end{equation*}
and, in fact, $\State_{n_1,n_2}(z)\doteq\hw{x,y}$ with $(x,y)$ given
by any pair from the Weyl orbit:
\begin{multline*}
  \bigl(\sqrt{p}-\ffrac{n_1}{\sqrt{p}} ,
  \sqrt{p}-\ffrac{n_2}{\sqrt{p}}\bigr),\
  \bigl(\ffrac{n_2}{\sqrt{p}}-\sqrt{p} ,
  \ffrac{n_1}{\sqrt{p}}-\sqrt{p}\bigr),
  \\
  \bigl(
  \sqrt{p}-\ffrac{n_2}{\sqrt{p}},
  \ffrac{n_1 + n_2}{\sqrt{p}}-2\sqrt{p}\bigr),\
  \bigl(\ffrac{n_1+n_2}{\sqrt{p}}-2 \sqrt{p},
  \sqrt{p}-\ffrac{n_1}{\sqrt{p}}\bigr),
  \\
  \bigl(\ffrac{n_1}{\sqrt{p}}-\sqrt{p},
  -\ffrac{n_1+n_2}{\sqrt{p}}+2 \sqrt{p}\bigr),\
  \bigl(-\ffrac{n_1+n_2}{\sqrt{p}}+2 \sqrt{p},
  \ffrac{n_2}{\sqrt{p}}-\sqrt{p}\bigr).
\end{multline*}
The corresponding $\algW_3$ singular vectors vanish 
in the free-field realization.  We propose the irreducible
$\algO$-modules generated from $\State_{n_1,n_2}(z)$ as counterparts
of the corresponding simple Yetter--Drinfeld $\Nich(X)$ modules, as a
starting point to study the relation between the two representation
categories.

\section{Conclusions}
We have outlined some details of the construction of the octuplet
extended algebra $\algO$ proposed in~\cite{[c-charge]}, and described
the corresponding Nichols algebra $\Nich(X)$ in rather explicit terms.
Systematically comparing $\algO$ representations with Yetter--Drinfeld
$\Nich(X)$ modules is very interesting from the perspective of whether
the relation existing in the $W_{p,1}$ (triplet-algebra)
case~\cite{[FGST],[FGST2],[NT],[TW]} extends to the current
$\algW_3$-related octuplet setting.

\subsubsection*{Acknowledgments} It is a pleasure to thank I.~Runkel
and C.~Schweigert for the useful discussions and the Department of
Mathematics, Hamburg University for hospitality.  This paper was
supported in part by the RFBR grant 11-01-00830.

\end{document}